\documentclass{article}

\usepackage{amsmath, amssymb, amsthm}
\usepackage{comment}
\usepackage{setspace}
\usepackage[backref]{hyperref}
\hypersetup{colorlinks=true}
\usepackage[colorinlistoftodos]{todonotes}
\usepackage{tikz}
\usetikzlibrary{arrows.meta}

	\newcommand{\setSize}[1]{\left|#1\right|}
		\newcommand{\Z}{\mathbb{Z}}
		
		\newcommand{\C}{\mathbb{C}}

		\newcommand{\pair}[1]{\left\langle #1 \right\rangle}
		\newcommand{\Laurent}[1]{(\!(#1)\!)}
		\newcommand{\Power}[1]{[\![#1]\!]}

	\newcommand{\indexSet}{I}
	\newcommand{\classicalIndex}{\indexSet_0}

    \newcommand{\affineGrassmannian}{\mathrm{Gr}}
    
    \newcommand{\maxParahoric}{G(\powerSeriesRing)}
    \newcommand{\below}[1]{B_{#1}}
	\newcommand{\iwahori}{\mathcal{I}}
	
	\newcommand{\powerSeriesRing}{\mathcal{O}}
	\newcommand{\laurentSeriesRing}{\mathcal{K}}

	\newcommand{\classicalRoots}{\roots_0}
	\newcommand{\positiveClassicalRoots}{\classicalRoots^+}
	
	\newcommand{\classicalSimpleRoots}{\simpleRoots_{0}}
	
	\newcommand{\classicalWeyl}{\weyl_{0}}
	
	\newcommand{\length}{\ell}

	\newcommand{\weyl}{W}
	\newcommand{\affineGrassmannians}{\weyl^{0}}
	\newcommand{\quantumBruhatGraph}{\Gamma}

	\newcommand{\rootLattice}{Q}
	\newcommand{\corootLattice}{\rootLattice^\vee}

		\newcommand{\roots}{\Phi}
		
		\newcommand{\simpleRoots}{\Delta}

	\newcommand{\mobius}{\mu}
	\newcommand{\relMobius}{\widetilde{\mobius}}
	\newcommand{\antidominantElements}{\tilde{\rootLattice}}

		\newcommand{\pathWeight}{\operatorname{wt}}
		\newcommand{\minPathWeight}{\operatorname{M}}

	\newcommand{\LieGroup}{G}

	\newcommand{\structureSheaf}[1]{O_{#1}}
	\newcommand{\idealSheaf}[1]{I_{#1}}

	\newcommand{\schubertCell}[1]{C_{#1}}
	\newcommand{\schubertVariety}[1]{X_{#1}}

	\newtheorem{lemma}{Lemma}
	\newtheorem{proposition}[lemma]{Proposition}
	\newtheorem{theorem}[lemma]{Theorem}
	\newtheorem{corollary}[lemma]{Corollary}
    \theoremstyle{remark}
    \newtheorem{example}[lemma]{Example}
    \newtheorem{remark}[lemma]{Remark}

\begin{document}

\title{The M\"obius Function on Affine Grassmannian elements}
\author{Michael Lugo and Mark Shimozono}

\maketitle

\begin{abstract}
    To any saturated chain in the affine Weyl group whose translation parts are sufficiently regular, we associate a near path and a far path in the quantum Bruhat graph. Using this,
    working in the Bruhat order on the minimal-length representatives of the cosets in the affine Weyl group with respect to the finite Weyl group, we characterize the pairs of elements for which the M\"obius function is nonzero. This is applied to obtain
    explicit expansions in the $K$-theory of affine Grassmannians, of the basis of ideal sheaves into the basis of structure sheaves of Schubert varieties.
\end{abstract}

\section{Introduction}

\subsection{M\"obius Function on Affine Grassmannian elements}
Let $\weyl$ be the Weyl group of a root system of untwisted affine type or the dual of untwisted affine type. For ease of exposition we shall work with untwisted affine type, leaving the statements of the dual case to \S \ref{S:dual}.
The group $\weyl$ is generated by simple reflections $s_i$ for $i$ in the affine Dynkin node set $\indexSet$, with distinguished affine node $0\in \indexSet$. Say that $w\in \weyl$ is \emph{affine Grassmannian} (denoted $w\in \affineGrassmannians$) if $ws_i>w$ for all $i\in \classicalIndex = \indexSet \setminus\{0\}$.
Let $\relMobius$ be the M\"obius function of the Bruhat order on $\affineGrassmannians$. The following is the affine Grassmannian case of Deodhar's general formula  \cite{deodhar} for the M\"obius function of a parabolic quotient of the Bruhat order. For $u,v\in \weyl$ with $u\le v$, denote by $[u,v]=\{z\in\weyl\mid u\le z\le v\}$ the interval between $u$ and $v$ in the Bruhat order on $\weyl$. 

\begin{theorem}[{\cite[Theorem 1.2]{deodhar}}] \label{T:deodhar}  For $u,v\in\affineGrassmannians$ with $u\le v$,
\begin{align}\label{deodhar_mobius}
\relMobius(u,v) = 
\begin{cases}
(-1)^{\ell(v)-\ell(u)} & \text{if $[u,v]\subset \affineGrassmannians$} \\
0 & \text{otherwise.}
\end{cases}
\end{align}
Moreover the first case holds if and only if there is no $i\in \classicalIndex$ such that $us_i\le v$. 
\end{theorem}

For $y\in \affineGrassmannians$ let $\below{y} = \{x\in \affineGrassmannians \mid [x,y]\subset\affineGrassmannians\}$.
Our main theorem (Theorem \ref{thm_explicit_result}) 
is an explicit description of the subposet $\below{y}$ of $\affineGrassmannians$ in the case that $y\in\affineGrassmannians$ is \emph{superregular}, that is, its translation part is sufficiently regular. For superregular $y$, we show that $\below{y}$  bijects with the classical Weyl group $\classicalWeyl$. Moreover, all of the possible poset structures on $\classicalWeyl$ which arise this way, come from a single directed graph structure $\quantumBruhatGraph$ on $\classicalWeyl$. This directed graph, first studied by Brenti, Fomin, and Postnikov \cite{bfp}, has become known as the quantum Bruhat graph 
due to its connection with quantum cohomology of flag varieties; it indicates which Schubert classes appear in the quantum product of a Schubert class with a Schubert divisor.

\subsection{Quantum Bruhat Graph}
For a positive classical root $\alpha\in \positiveClassicalRoots$ let $r_\alpha\in \classicalWeyl$ be the associated reflection and $\alpha^\vee$ the associated coroot. Let $\rho=(1/2) \sum_{\alpha\in \positiveClassicalRoots} \alpha$.

The quantum Bruhat graph (QBG) \cite{bfp} is the directed graph $\quantumBruhatGraph$ with vertex set $W_0$ and a directed edge $w\to wr_\alpha$ for each pair $(w,\alpha)\in \classicalWeyl\times  \positiveClassicalRoots$ such that one of the following holds.
\begin{enumerate}
	\item $\length(wr_\alpha) = \length(w)+1$. Such edges are called \emph{Bruhat} edges.
	\item $\length(wr_\alpha) = \length(w) - \pair{\alpha^\vee,2\rho}+1$. Such edges are called \emph{quantum} edges.
\end{enumerate}
The \emph{label} of the edge $w\to wr_\alpha$ is $\alpha\in \positiveClassicalRoots$.
The \emph{weight} of the edge $w\to wr_\alpha$ in $\quantumBruhatGraph$ is $\alpha^\vee$ if the edge is quantum and $0$ if it is Bruhat.
The \emph{weight} of a path\footnote{All our paths are directed.} $P$, denoted $\pathWeight(P)$, is the sum of the weights of its edges. 
For any $u, v \in \classicalWeyl$, let $\minPathWeight(u,v)$ be the weight of any shortest path in $\quantumBruhatGraph$ from $u$ to $v$; this is well defined by \cite[Lemma 1]{postnikov}.

\begin{example} The quantum Bruhat graph of type $A_2$ is pictured in Figure \ref{F:QBG}. Let $w_0=s_1s_2s_1$ and $\alpha_{12}^\vee=\alpha_1^\vee+\alpha_2^\vee$. 
The quantum edges are given in red and their weights are indicated. We have $\minPathWeight(w_0,s_1)=\alpha_1^\vee+\alpha_2^\vee$, which is realized by either of the shortest paths 
$w_0 \to \mathrm{id} \to s_1$ or $w_0\to s_1s_2 \to s_1$.
\begin{figure}
\[
\begin{tikzpicture}
\draw (0:0) node (w0) {$w_0$} ++ (330:2) node (s21) {$s_2s_1$} ++ (270:2) node (s2) {$s_2$} ++ (210:2) node (id) {$\mathrm{id}$}
++ (150:2) node (s1) {$s_1$} ++ (90:2) node (s12) {$s_1s_2$};
\path[->] (id) edge [bend right=15] (s1);
\path[->,color=red] (s1) edge [bend right=15] node [below,color=black]  {$\alpha_1^\vee$} (id);
\path[->] (id) edge [bend left=15] (s2);
\path[->,color=red] (s2) edge [bend left=15] node[below,color=black] {$\alpha_2^\vee$} (id);
\draw[->] (s1) edge [bend right=15]  (s12);
\path[->,color=red] (s12) edge [bend right=15] node[left,color=black] {$\alpha_2^\vee$}(s1);
\draw[->] (s1) -- (s21);
\draw[->](s2) edge [bend left=15] (s21);
\draw[->] (s2) -- (s12);
\draw[->,color=red] (s21) edge [bend left=15] node [right,color=black] {$\alpha_1^\vee$} (s2);
\draw[->] (s12)  edge [bend right=15] (w0);
\draw[->](s21) edge [bend left=15] (w0);
\draw[->,red] (w0)--(id) node [pos=.25,right=-2pt,color=black] {$\alpha_{12}^\vee$};
\draw[->,red] (w0) edge [bend right=15] node [above,color=black] {$\alpha_1^\vee$} (s12);
\draw[->,red] (w0) edge [bend left=15] node [above,color=black] {$\alpha_2^\vee$} (s21);
\end{tikzpicture}
\]
\caption{Quantum Bruhat Graph of type $A_2$}
\label{F:QBG}
\end{figure}
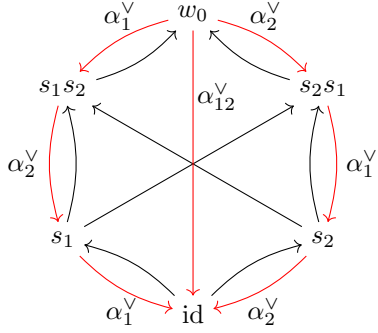
\end{example}

\subsection{Main Theorem}
An element $\lambda$ of the classical coroot lattice $\corootLattice$ is \emph{superregular} if $|\pair{\lambda,\alpha_i}| \gg0$\footnote{Lower bounds for superregularity are discussed in \S \ref{sec_regularity}.} for all $i\in \classicalIndex$.
Any element $x\in \weyl$ can be written uniquely in the form $x=wt_\lambda$ with $w\in\classicalWeyl$ and $\lambda\in\corootLattice$ where
$t_\lambda\in\weyl$ is the translation element. For $w\in \classicalWeyl$ we define the affine Weyl group element $wt_\lambda$ to be superregular if $\lambda$ is.

\begin{theorem}\label{thm_explicit_result}
	Let $x,y\in \affineGrassmannians$ with $x = w't_{\lambda'}$, $y=wt_{\lambda}$ with $w,w'\in \classicalWeyl$, $\lambda,\lambda'\in\corootLattice$ and $\lambda$ superregular. Then
	\[
	\relMobius(x,y) = \left\{
	\begin{array}{ll}
	(-1)^{\length(y)-\length(x)} &\text{if } \lambda' = \lambda + \minPathWeight(w,w') \\
	0 &\text{else}
	\end{array}
	\right.
	\]
	where $\minPathWeight$ is the weight of a shortest path in $\quantumBruhatGraph$ from $w$ to $w'$.
\end{theorem}

\begin{corollary} Let $y=wt_\lambda\in \affineGrassmannians$ be superregular. Then the set
	$$\below{y} = \{ x\in \affineGrassmannians \mid \text{$x\le y$ and $\relMobius(x,y)\ne0$} \}$$ is given by 
\begin{align}\label{E:elements below}
	\{ u t_{\lambda+\minPathWeight(w,u)} \mid u\in \classicalWeyl\}.
\end{align}
\end{corollary}

\begin{remark}
The poset structure on $\below{y}$ can be constructed from the shortest paths from $w$ to all the other 
vertices in $\quantumBruhatGraph$. This poset structure is given by the 
$w$-tilted Bruhat order $\preceq_w$ of \cite{bfp}, which is defined by
$u\preceq_w v$ if $u$ is on some shortest path in $\quantumBruhatGraph$ from $w$ to $v$.
\end{remark}

\subsection{Application to the \texorpdfstring{$K$}{K}-theory of the affine Grassmannian}
Let $\powerSeriesRing=\C\Power{t}$ be the ring of formal power series with coefficients in $\C$,
$\laurentSeriesRing=\C\Laurent{t}$ the field of formal Laurent series,
$\LieGroup$ a simple Lie group over $\C$, $\affineGrassmannian_G =G(\laurentSeriesRing)/G(\powerSeriesRing)$ the affine Grassmannian, $\iwahori$ the Iwahori subgroup of the affine Kac-Moody group. For $x\in \affineGrassmannians$ let  $\schubertCell{x}=\iwahori x \maxParahoric/\maxParahoric\subset \affineGrassmannian_G$ be the affine Grassmannian Schubert cell and $\schubertVariety{x} = \overline{C_x}\subset \affineGrassmannian_G$ the Schubert variety.
For $x\in \affineGrassmannians$ let $\structureSheaf{x}$ and $\idealSheaf{x}$ be the elements of $K(\affineGrassmannian_G)$ given by the $K$-classes of the structure sheaf 
of $\schubertVariety{x}$ and the ideal sheaf of the boundary $\partial{\schubertVariety{x}}$.

The following is due to Kumar \cite{kumar_email}, who produced the proof (for Kac-Moody partial flag varieties) upon being asked about it by the second author.

\begin{proposition}[\cite{kumar_email}] For any $y\in\affineGrassmannians$
	\begin{align}
	O_y = \sum_{\substack{x\in \affineGrassmannians \\ x \le y }} I_x.
\end{align}
\end{proposition}

By the definition of the M\"obius function and Theorem \ref{thm_explicit_result} we deduce:

\begin{corollary}
	Let $y = wt_\lambda \in \affineGrassmannians$ be superregular. Then
	\[
	I_y = \sum_{u \in \classicalWeyl}(-1)^{\length(w,u)} O_{ut_{\lambda + \minPathWeight(w,u)}}
	\]
	where $\length(w,u)$ is the length of a shortest path from $w$ to $u$ in $\quantumBruhatGraph$.
\end{corollary}

\section{Near and far paths in \texorpdfstring{$\quantumBruhatGraph$}{the QBG} associated to superregular saturated chains in \texorpdfstring{$\weyl$}{\weyl}}

\subsection{Sufficiently Regular Bruhat covers}
	
	Say that $\lambda \in \corootLattice$ is \emph{sufficiently regular} if $|\pair{\lambda, \alpha_i}| \gg 0$ for all $i \in \classicalIndex$ and that $y = wt_\lambda\in\weyl$ is sufficiently regular if $\lambda$ is\footnote{This definition follows the definition of superregular given in \cite{lam_shimo}. In this paper, superregularity has a different requirement.}. The question of sufficient regularity to qualify will be addressed in \S \ref{sec_regularity}. Say that $\lambda\in \corootLattice$ is \emph{antidominant} if $\pair{\lambda,\alpha_i} \le 0$ for all $i\in \classicalIndex$.
	Let $\antidominantElements$ be the set of antidominant elements in $Q^\vee$.

	\begin{proposition}[{\cite[Proposition 4.1]{lam_shimo}}]\label{prop_the_big_four}
		Let $\lambda \in \antidominantElements$ be sufficiently regular and $y = wt_{v\lambda}$ where $w,v\in\classicalWeyl$. Then $x = yr_{v\alpha + n\delta} \lessdot y$ for some $\alpha\in \positiveClassicalRoots$ and $n\in\Z$ if and only if one of the following conditions hold:
		\begin{enumerate}
			\item $\length(wv) = \length(wvr_\alpha) - 1$ and $n = \pair{\lambda, \alpha}$, giving $x = wr_{v\alpha}t_{v\lambda}$.
			\item $\length(wv) = \length(wvr_\alpha) + \pair{\alpha^\vee, 2\rho} - 1$ and $n = \pair{\lambda, \alpha} + 1$, giving $x = wr_{v\alpha}t_{v(\lambda + \alpha^\vee)}$.
			\item $\length(v) = \length(vr_\alpha) + 1$ and $n=0$, giving $x = wr_{v\alpha}t_{vr_\alpha(\lambda)}$
			\item $\length(v) = \length(vr_\alpha) - \pair{\alpha^\vee, 2\rho} + 1$ and $n = -1$, giving $x = wr_{v\alpha} t_{vr_{\alpha}(\lambda + \alpha^\vee)}$.
		\end{enumerate}
		\end{proposition}
	That is, every sufficiently regular cover in the Bruhat order on $\weyl$ has a corresponding edge in $\quantumBruhatGraph$.
	In the terminology introduced in \cite{lam_shimo}, the first two cases are called \emph{near}, since the chamber of the translation component remains fixed, and the last two cases \emph{far}.

	\begin{corollary}\label{cor_the_two} The four possibilities in the previous Proposition may be replaced by the following two:
		\begin{enumerate}
			\item $wv \to wvr_\alpha$ is an edge $E$ in $\quantumBruhatGraph$ and $x = wr_{v\alpha}t_{v(\lambda + \pathWeight(E))}$. We call this a \emph{near covering}.
			\item $vr_\alpha \to v$ is an edge $E$ in $\quantumBruhatGraph$ and $x = wr_{v\alpha} t_{vr_{\alpha}(\lambda + \pathWeight(E))}$. We call this a \emph{far covering}.
		\end{enumerate}
		\end{corollary}

\subsection{Superregular saturated Bruhat chains}
	We extend this idea to saturated chains in the Bruhat order. We say that $y = wt_\lambda \in \weyl$ is \emph{superregular} with respect to $x$ if for every saturated chain $x \lessdot z_{k-1} \lessdot \dots \lessdot z_1 \lessdot y$ it is true that $y$ and $z_i$ for $i = 1,\dots, k-1$ are all sufficiently regular. We say a saturated chain $C$ is superregular if its top element is superregular with respect to its bottom element. Suppose that 
	\begin{align}\label{E:saturated chain}
	x \xrightarrow{\beta_k} z_{k-1} \xrightarrow{\beta_{k-1}} \dots \xrightarrow{\beta_2} z_1 \xrightarrow{\beta_1} y
	\end{align}
	is a superregular saturated chain in $\weyl$, where $\beta_i \in \positiveClassicalRoots$ is the root corresponding to $\alpha$ in Proposition \ref{prop_the_big_four} which labels the corresponding edge in $\quantumBruhatGraph$. Let $1 \leq f_1<f_2<\dots<f_i \leq k$ be the ordered subindices relating to far coverings and $1 \leq n_1<n_2<\dots<n_j \leq k$ be the ordered subindices relating to near coverings. 
	Define the \emph{far product} $r_f = r_{\beta_{f_1}}r_{\beta_{f_2}}\dots r_{\beta_{f_i}}$ and \emph{near product} $r_n = r_{\beta_{n_1}}r_{\beta_{n_2}}\dots r_{\beta_{n_j}}$.

	\begin{lemma}\label{lem_near_far_wherever_you_are}
		Suppose $x < y=wt_{v\lambda}$ with $\lambda\in\antidominantElements$ and superregular chain as in Equation \eqref{E:saturated chain}. Let $r_f$ and $r_n$ be defined as above. Then there exist well defined paths $P_f: vr_f \to v$ and $P_n: wv \to wvr_n$ such that
		\[
			x =
			wvr_n(vr_f)^{-1}t_{vr_f\left(\lambda + \pathWeight(P_n) + \pathWeight(P_f)\right)}
		\]
		\end{lemma}
		\begin{example}\label{example_chain}
		    Let $1$ denote the identity element of $\classicalWeyl$. Consider the following superregular saturated chain in type $A_2^{(1)}$.
		    \begin{align*}
		        x &= s_2t_{s_1s_2s_1(-2\alpha_1^\vee-3\alpha_2^\vee)}
		        \xrightarrow{\alpha_1+\alpha_2} \\
		        z_2 &=s_1s_2t_{-3\alpha_1^\vee-4\alpha_2^\vee}
		        \xrightarrow{\alpha_1} \\
		        z_1 &= s_1s_2s_1t_{-4\alpha_1^\vee-4\alpha_2^\vee}
		        \xrightarrow{\alpha_1} \\
		        y &= s_1s_2t_{-4\alpha_1^\vee-4\alpha_2^\vee}
		    \end{align*}
		    We observe that $x \lessdot z_2$ is a far covering that bijects with the quantum edge $s_1s_2s_1 \to 1$ with path weight $\alpha_1^\vee + \alpha_2^\vee$. $z_2 \lessdot z_1$ is a near covering that bijects with the quantum edge $s_1s_2s_1 \to s_1s_2$ with path weight $\alpha_1^\vee$. $z_1 \lessdot y$ is a near covering that bijects with the Bruhat edge $s_1s_2 \to s_1s_2s_1$ with path weight 0.
		    
		    In terms of Lemma \ref{lem_near_far_wherever_you_are}, this implies that
		    \begin{itemize}
		        \item $w = s_1s_2$, $v = 1$, and $\lambda = -4\alpha_1^\vee-4\alpha_2^\vee$
		        \item $r_f = s_1s_2s_1$ and $r_n = (s_1)(s_1) = 1$
		        \item $P_f: s_1s_2s_1 \to 1$ with $\pathWeight(P_f) = \alpha_1^\vee+\alpha_2^\vee$
		        \item $P_n: s_1s_2 \to s_1s_2$ with $\pathWeight(P_n) = \alpha_1^\vee$
		    \end{itemize}
		    
		    We find $x$ can be expressed in these terms.
		    \begin{align*}
		        x 
		        &=wvr_n(vr_f)^{-1}t_{vr_f(\lambda + \pathWeight(P_n) + \pathWeight(P_f))} \\
		        &= (s_1s_2)(1)(1)[(1)(s_1s_2s_1)]^{-1}t_{(1)(s_1s_2s_1)(-4\alpha_1^\vee-4\alpha_2^\vee + \alpha_1^\vee + \alpha_1^\vee + \alpha_2^\vee)} \\
		        &= s_1s_2s_1s_2s_1t_{s_1s_2s_1(-2\alpha_1^\vee-3\alpha_2^\vee)} \\
		        &= s_2t_{s_1s_2s_1(-2\alpha_1^\vee-3\alpha_2^\vee)}
		    \end{align*}
		\end{example}
		
		\begin{proof}
			We proceed by induction. The base case for $k = 1$ is given explicitly by Corollary \ref{cor_the_two}: the near case yields $P_f$ as trivial while $P_n$ is the edge described, and the far case yields $P_n$ as trivial while $P_f$ is the edge described. Suppose the result is true for any fixed length difference $k\geq1$ and that $\length(y) - \length(x) = k+1$. By induction, there are elements $r'_n, r'_f$ and associated paths $P'_n: wv \to wvr'_n, P'_f: vr'_f \to v$ such that $z_k = wvr'_n(vr'_f)^{-1}t_{vr'_f(\lambda + \pathWeight(P_n') + \pathWeight(P'_f))}$. Let $\beta = \beta_{k+1}$. We will perform cases based on whether the covering $x \lessdot z_k$ is a near or far covering. Note that in either case, we observe the following is the resulting classical Weyl component of $x$:
			\[
				wvr_n'(vr_f')^{-1}r_{vr_f'(\beta)}
				= wvr_n'(vr_f')^{-1}(vr_f')r_\beta(vr_f')^{-1}
				= wvr_n'r_\beta(vr_f')^{-1}
			\]
			If $x \lessdot z_k$ by a near covering, then $r'_f = r_f$, $r'_nr_{\beta} = r_n$,  and $P'_f = P_f$. The covering yields an edge $E$ from $wvr'_n(vr'_f)^{-1}(vr'_f) = wvr'_n \to wvr'_nr_\beta = wvr_n$ and shows that $P'_n + E = P_n$ and that $\pathWeight(P'_n)+\pathWeight(E) = \pathWeight(P_n)$. Therefore
			\begin{align*}
				x
				&= wvr_n'r_\beta(vr_f')^{-1} t_{vr'_f(\lambda + \pathWeight(P_n') + \pathWeight(P_f') + \pathWeight(E))} \\
				&= wvr_n'r_\beta(vr_f)^{-1}t_{vr_f(\lambda + \pathWeight(P_n') + \pathWeight(P_f) + \pathWeight(E))} \\
				&= wvr_n(vr_f)^{-1}t_{vr_f(\lambda + \pathWeight(P_n) + \pathWeight(P_f))}
			\end{align*}
			If $x \lessdot z_k$ by a far edge, then $r'_n = r_n$, $r'_fr_\beta=r_f$, and $P'_n = P_n$. The covering yields an edge $E$ from $vr'_fr_\beta = vr_f \to vr'_f$ and shows that $E + P'_f = P_f$ and $\pathWeight(P'_f) + \pathWeight(E) = \pathWeight(P_f)$. Therefore
			\begin{align*}
				x
				&=wvr_n'r_\beta(vr_f')^{-1}t_{vr'_fr_{\beta}(\lambda + \pathWeight(P_n') + \pathWeight(P_f') + \pathWeight(E))} \\
				&=wvr_nr_\beta {r_f'}^{-1}v^{-1}t_{vr'_fr_{\beta}(\lambda + \pathWeight(P_n) + \pathWeight(P_f') + \pathWeight(E))} \\
				&=wvr_nr_f^{-1}v^{-1}t_{vr_f(\lambda + \pathWeight(P_n) + \pathWeight(P_f))} \\
				&=wvr_n(vr_f)^{-1}t_{vr_f(\lambda + \pathWeight(P_n) + \pathWeight(P_f))}
			\end{align*}
			So by induction, the result is true for all superregular saturated chains. 
			\end{proof}

        For a superregular saturated chain \eqref{E:saturated chain}, we call $P_f$ the associated \emph{far path} and $P_n$ the associated \emph{near path}.

    The following Corollary is the specialization of 
    Lemma \ref{lem_near_far_wherever_you_are} to $\affineGrassmannians$.

	\begin{corollary}\label{cor_near_far_specialized}
		Suppose $x,y \in \affineGrassmannians$, $x < y=wt_{\lambda}$ and let \eqref{E:saturated chain} be a superregular saturated chain where $\beta_i$, $r_f, r_n, P_f,P_n$ are defined as in Lemma \ref{lem_near_far_wherever_you_are}. Then $P_f$ goes from $1$ to $1$, $P_n$ goes from $w$ to $wr_n$, and
		\[
			x =
			wr_nt_{\lambda + \pathWeight(P_n) + \pathWeight(P_f)}
		\]
		\end{corollary}
		\begin{proof}
			Since $y \in \affineGrassmannians$, $v=1$ in Lemma \ref{lem_near_far_wherever_you_are}. Similarly, since $x \in \affineGrassmannians$, $v r_f = r_f = 1$. Therefore $P_f$ is a closed loop at the identity. The proof is complete.
			\end{proof}
			
	\subsection{Loops of length 2 in \texorpdfstring{$\quantumBruhatGraph$}{the QBG}}
	    We note two facts about the Bruhat order.
	
    	\begin{theorem}[{\cite[Chain Property]{comb_of_cox}}]\label{thm_chain_property}
    		If $u < w$ and $u,w\in \affineGrassmannians$, then there exists a saturated chain 
    		$u=w_0\lessdot w_1 \lessdot w_2 \lessdot \dotsm \lessdot w_k = w$ with $w_i\in \affineGrassmannians$.
    	\end{theorem}
    	
    	\begin{lemma}[{\cite[Lemma 2.7.3]{comb_of_cox}}]\label{lem_small_interval}
    		If $x<y$ in $\weyl$ with $\length(y) = \length(x) + 2$, then $[x,y] = \{x,u,z,y\}$ consists of four distinct elements with $x \lessdot u \lessdot y$ and $x \lessdot z \lessdot y$.
    	\end{lemma}

        We observe that any loop of length 2 in $\quantumBruhatGraph$ has a simple root as the label for both of its edges; apply both quantum and Bruhat edge conditions to $\alpha$.

    	\begin{lemma}\label{lem_bad_loops}
    		Let $x, y \in \affineGrassmannians$ with $y$ superregular, $\length(y) = \length(x) + 2$ and $x < y$. 
    		Let $z\in\affineGrassmannians$ be such that $x\lessdot z\lessdot y$ is the chain guaranteed to exist by Theorem \ref{thm_chain_property}. Then $[x,y] \not \subset \affineGrassmannians$ if and only if the near path associated to $x \lessdot z \lessdot y$ is a 2-loop.
		\end{lemma}
		\begin{proof}
			Let $y = wt_\lambda$. Suppose $u \not\in\affineGrassmannians$. Then by Corollary \ref{cor_near_far_specialized}, $x \lessdot u \lessdot y$ associates to a trivial near path, a far path that is a loop of length 2 at the identity with weight $\alpha^\vee$ for some $\alpha \in \classicalSimpleRoots$, $u= ws_\alpha t_{s_\alpha(\lambda)}$, and $x = wt_{\lambda+\alpha^\vee}$. Thus by Corollary \ref{cor_near_far_specialized}, the associated near path for $x \lessdot z \lessdot y$ is a loop of length 2 at $w$. Now suppose the near path associated to $x \lessdot z \lessdot y$ is a loop $P_n: w \to ws_\alpha \to w$ ($\alpha \in \classicalSimpleRoots$). Then $x = wt_{\lambda+\alpha^\vee}$ by Corollary \ref{cor_near_far_specialized}. If $u \in \affineGrassmannians$, this and Lemma \ref{cor_near_far_specialized} would imply a different $\beta^\vee$ such that
			\[
				wt_{\lambda + \beta^\vee}
				=x 
				= wt_{\lambda + \alpha^\vee}
				\Rightarrow \alpha= \beta
			\]
			But this contradicts $u \not = z$.
		\end{proof}
		
		\begin{example}
		    In Example \ref{example_chain}, we find that $z_1 \lessdot z_2 \lessdot y$ corresponds with a near path that is a 2-loop at $s_1s_2$ and $z_1,z_2,y \in \affineGrassmannians$. In this case, $u = s_1s_2s_1t_{s_1(-3\alpha_1^\vee-4\alpha_2^\vee)}$ satisfies $u \not\in\affineGrassmannians$ and $u \in [z_1,y]$.
		\end{example}
		
		\begin{remark}\label{remark_exchanges}
		    It is always the case that if the near path 2-loop from Lemma \ref{lem_bad_loops} utilizes the simple root $\alpha$, then the associated far path which induces $u$ is the far path is the loop $1 \to r_\alpha \to 1$.
		\end{remark}

    \subsection{Interval equivalent paths}
        
        We now discuss and extend properties of paths in $\quantumBruhatGraph$ uncovered by Postnikov \cite{postnikov}.
        Following Dyer \cite{dyer}, a reflection ordering is equivalent to a total order on $\positiveClassicalRoots$ such that for all $\alpha, \beta \in \classicalSimpleRoots$ with $\alpha < \beta$ and $\alpha + \beta \in \positiveClassicalRoots$, then $\alpha < \alpha + \beta < \beta$.
        Given a reflection ordering on $\positiveClassicalRoots$, then the following is true.

        \begin{lemma}[Lemma 6.7 \cite{bfp}]\label{lem_diamond}
            If $a \xrightarrow{\alpha_1} x \xrightarrow{\beta_1} c$ is a path in $\quantumBruhatGraph$ with $\alpha_1 > \beta_1$, then there exists a unique $y \in \classicalWeyl$ such that $a \xrightarrow{\alpha_2} y \xrightarrow{\beta_2} c$ where $\beta_1 < \beta_2 > \alpha_2 < \alpha_1$.
        \end{lemma}
        
        
        In essence, any descent in the label sequence of a path can be uniquely exchanged for an ascent without changing the endpoints of the path. Note that the reversed ordering of a reflection ordering ($\alpha < \beta$ in the reversed ordering if and only if $\beta < \alpha$ in the original reflection ordering) is itself a reflection ordering. By applying Lemma \ref{lem_diamond} to the reversed ordering, we find that we may also uniquely exchange any descent in the label sequence of a path with an ascent.
        
        Given a path $P$, we can apply these swaps to sort the label sequence of a path. Pairing this information with \cite[Theorem 6.4]{bfp}, we obtain two key facts.
        
        \begin{corollary}[{\cite[proof of Lemma 1]{postnikov}, \cite{bfp}}]\label{cor_loops_appear}
            Suppose $P:w \to w'$ is a non-minimal path in $\quantumBruhatGraph$. 
        	\begin{enumerate}
        	    \item $P$ is interval equivalent to a path containing a loop of length 2.
        	    \item $\pathWeight(P) - \minPathWeight(w, w')$ is a positive sum of coroots.
        	\end{enumerate}
        \end{corollary}
        
        We say that two paths in $\quantumBruhatGraph$ are \emph{interval equivalent} if the paths can be obtained from each other by a series of exchanges using Lemma \ref{lem_diamond}.
        In Postnikov's proof of \cite[Lemma 1]{postnikov}, he notes Lemma \ref{lem_diamond} also preserves the weight of the path.

        \begin{lemma}\label{lem_interval_equivalence}
            Let $C$ be a superregular saturated chain from $x$ to $y$ which utilizes only near (resp. far) coverings, $P$ be the associated near (resp. far) path, and $P'$ be interval equivalent to $P$. Then there is a unique superregular saturated chain $C'$ from $x$ to $y$ which utilizes only near (resp. far) coverings with corresponding near (resp. far) path $P'$.
        \end{lemma}
        \begin{proof}
            By Proposition \ref{prop_the_big_four}, each edge of $P$ bijects to a covering in $C$. By Lemma \ref{lem_near_far_wherever_you_are}, the exact element $x$ can be found using the weight and endpoint of $P$. If $P'$ is interval equivalent, then $P'$ shares the same endpoint, length, and weight of $P$. Thus, if we biject the edges of $P'$ with near (resp. far) coverings then Lemma \ref{lem_near_far_wherever_you_are} states that the lower element from the resulting saturated chain $C'$ must be $x$ as well.
        \end{proof}
        
        This previous lemma is sufficient to prove Theorem $\ref{thm_explicit_result}$. However, it can be strengthened by the following proposition.
        
        \begin{proposition}
            Let $C$ be a superregular saturated chain from $x$ to $y$, $P_n$ be the associated near path, and $P_f$ be the associated far path. If $P_n$ is interval equivalent to $P_n'$, and $P_f$ is interval equivalent to $P_f'$, then
            \begin{enumerate}
                \item there exists a saturated chain $C'$ from $x$ to $y$ with associated near and far paths $P_n'$ and $P_f'$.
                \item any saturated chain $C'$ with upper element $y$, associated near path $P_n'$, and associated far path $P_f'$ has lower element $x$.
            \end{enumerate}
        \end{proposition}
        \begin{proof}
            Note that the proof by induction of Lemma \ref{lem_near_far_wherever_you_are} shows that each near covering only modifies $r_n$ and $P_n$ while leaving the far components fixed. Similarly, each far covering only modifies $r_f$ and $P_f$. Thus, having a near covering then far covering results in the same lower element as when we take the same far covering first and then the near covering. The only restriction is one cannot change the order of near coverings with respect to other near coverings and similarly for far coverings.
            
            To show (1) and (2), we may take the saturated chain $C_1$ given by Lemma \ref{lem_interval_equivalence} from $z$ to $y$ generated by $P_n'$, then append the saturated chain $C_2$ from $x$ to $z$ given by $P_f'$.
        \end{proof}
        
        \begin{example}
            Consider the sub-chain $x \xrightarrow{\alpha_1+\alpha_2} z_2 \xrightarrow{\alpha_1} z_1$ from Example \ref{example_chain}. The left covering was a far covering corresponding to the quantum edge $s_1s_2s_1 \to 1$ while the right covering was a near covering corresponding to the quantum edge $s_1s_2s_1 \to s_1s_2$. If instead we used the quantum edge $s_1s_2s_1 \to 1$ of a far covering with upper element $z_1$, we generate the element $z_3 = t_{s_1s_2s_1(-3\alpha_1^\vee-3\alpha_2^\vee)}$. We also observe that $x \xrightarrow{\alpha_1} z_3$ by a near covering. This demonstrates that changing the order of near and far covers with respect to each other results in a chain still contained in the interval: $z_3 \in [x,z_1]$.
        \end{example}
        
\subsection{Proof of Theorem \ref{thm_explicit_result}}

	Note that proving Theorem \ref{thm_explicit_result} is equivalent to showing that for $x=w't_{\lambda'},y=wt_\lambda \in \affineGrassmannians$ and $y$ superregular with respect to $x$, then $[x,y] \subset \affineGrassmannians$ if and only if $\lambda' = \lambda + \minPathWeight(w,w')$.
	
	Suppose $[x,y] \subset \affineGrassmannians$. Theorem \ref{thm_chain_property} implies there exists a saturated chain in $\affineGrassmannians$ from $x$ to $y$, which by Corollary \ref{cor_near_far_specialized} implies $x = w't_{\lambda + \pathWeight(P_n)}$ for the associated near path. If $P_n$ has non-minimal length, then by Corollary \ref{cor_loops_appear} part 1, there exists a $P_n'$ which is interval equivalent to $P_n$ but contains a loop of length 2. However, Lemma \ref{lem_bad_loops} implies the induced saturated chain in $[x,y]$ is not contained in $\affineGrassmannians$, which is a contradiction. Thus $P_n$ has minimal length, i.e. $\pathWeight(P_n) = \minPathWeight(w,w')$.
	
	Suppose $\lambda' = \lambda + \minPathWeight(w,w')$.
	If $[x,y] \not\subset\affineGrassmannians$, then there exists a saturated chain from $x$ to $y$ which utilizes far coverings. Corollary \ref{cor_near_far_specialized} supplies that $\lambda' = \lambda + \pathWeight(P_n) + \pathWeight(P_f)$ for the associated near and far paths. $P_f$ must be a non-empty loop at the identity and $P_n: w \to w'$. But this implies $\pathWeight(P_n) - \minPathWeight(w,w') = -\pathWeight(P_f)$ with $\pathWeight(P_f) \not = 0$ since $P_f$ is a non-empty loop, which is a contradiction to Corollary \ref{cor_loops_appear} part 2. \qed

\section{Bounds on Regularity}\label{sec_regularity}
	
	Initially, the definition of sufficiently regular in \cite{lam_shimo} left a large bound of $|\pair{\lambda, \alpha}| > 2|\classicalWeyl| + 2$ for all $\alpha \in \classicalSimpleRoots$. Later, in \cite[Proposition 4.2]{milicevic}, it was found that one only requires $|\pair{\lambda, \alpha}| \geq 2\length(w_0)$ if $\classicalWeyl \not= G_2$ (change to $3\length(w_0)$ in the case of $G_2$) for all such $\alpha \in \classicalSimpleRoots$. Recently in \cite{welch}, it was shown that Proposition \ref{prop_the_big_four} holds when $\classicalWeyl$ is simply laced and $|\pair{\lambda, \alpha}| \geq 3$ for all $\alpha \in \classicalSimpleRoots$. Furthermore, this hints at a similar, smaller bound for the non-simply laced cases using a similar proof to \cite{milicevic}. For the following discussion, suppose the regularity bound for Proposition \ref{prop_the_big_four} is $k$.
	
	We can say that $y$ is superregular with respect to $x$ if we are able to use Proposition \ref{prop_the_big_four} for every covering between $x$ and $y$ in a saturated chain. If $\length(y) - \length(x) = m$, then we must use Proposition \ref{prop_the_big_four} $(m-1)$-times. We note that each use of Proposition \ref{prop_the_big_four} can modify the translation by adding $\alpha^\vee$. Thus we can decrease the regularity by at most the maximum of $\pair{\beta, \alpha^\vee}$ over $\beta \in \classicalSimpleRoots$. If $j$ is the mentioned maximum, then one functioning bound for $\lambda$ is $k + (m-1)j$.
	
	Theorem \ref{thm_explicit_result} requires the use of Proposition \ref{prop_the_big_four} up to the length of the longest shortest path starting at $w$. Note that $\setSize{\classicalWeyl}$ is an upper bound on the longest shortest path. One functioning lower bound on the regularity of $\lambda$ in Theorem \ref{thm_explicit_result} is $k + \setSize{\classicalWeyl}j$.e
	

\section{Dual untwisted affine root systems}
\label{S:dual}
Consider the dual of an untwisted affine root system. There is an associated quantum Bruhat graph defined in the same way as for the untwisted affine root systems, except that the roles of roots and coroots are exchanged.
For dual untwisted type $s_0 = t_\varphi s_\varphi$ where $\varphi$ is the short dominant root.

Let $\rho^\vee=(1/2) \sum_{\alpha^\vee\in {\positiveClassicalRoots}^\vee} \alpha^\vee$.

The quantum Bruhat graph (QBG) is the directed graph $\quantumBruhatGraph$ with vertex set $\classicalWeyl$ and a directed edge $w\to wr_\alpha$ for each pair $(w,\alpha)\in \classicalWeyl\times  \positiveClassicalRoots$ such that one of the following holds.
\begin{enumerate}
	\item (Bruhat edge) $\length(wr_\alpha) = \length(w)+1$. 
	\item (Quantum edge) $\length(wr_\alpha) = \length(w) - \pair{2\rho^\vee,\alpha}+1$.
\end{enumerate}
The \emph{weight} of the edge $w\to wr_\alpha$ in $\quantumBruhatGraph$ is $\alpha$ if the edge is quantum and $0$ if it is Bruhat.
For any $u, v \in \classicalWeyl$, let $\minPathWeight(u,v)$ be the weight of any shortest path in $\quantumBruhatGraph$ from $u$ to $v$.

With these definitions in place, analogues of the statements for the untwisted affine root systems hold for the dual root systems with essentially the same proofs.

\begin{example} 
\begin{enumerate}
\item 
Consider the root system $A_{2\cdot3-1}^{(2)}$, which is dual to $B_3^{(1)}$.

\begin{figure}
\[
\tikzstyle{dynkin}=[circle,draw=black,thick,
                   inner sep=0pt,minimum size=5mm]
\begin{matrix}
\begin{tikzpicture}
\node[dynkin] (n0) at (0,1) {$0$};
\node[dynkin] (n1) at (0,-1) {$1$};
\node[dynkin] (n2) at (1,0) {$2$};
\node[dynkin] (n3) at (2,0) {$3$};
\draw[thick] (n0) -- (n2) -- (n1);
\draw[-Implies,line width=1pt,double distance=1pt] (n3) -- (n2);
\end{tikzpicture} &&
\begin{tikzpicture}
\node[dynkin] (n0) at (0,0) {$0$};
\node[dynkin] (n1) at (1,0) {$1$};
\node[dynkin] (n2) at (2,0) {$2$};
\node[dynkin] (n3) at (3,0) {$3$};
\node[] (nx) at (0,-1) {\vphantom{$\overline{0}$}};
\draw[-Implies,line width=1pt,double distance=1pt] (n1) -- (n0);
\draw[thick] (n1) -- (n2);
\draw[-Implies,line width=1pt,double distance=1pt] (n2) -- (n3);
\end{tikzpicture} \\
A_{2\cdot3-1}^{(2)} & &D_{3+1}^{(2)}
\end{matrix}
\]
\label{F:Aodd_twisted}
\caption{Some dual untwisted affine Dynkin diagrams}
\end{figure}
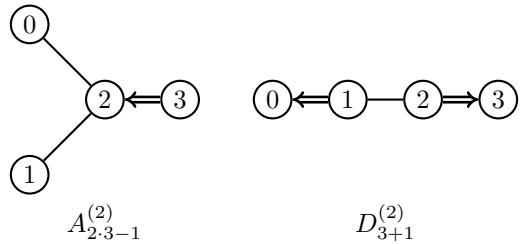
The classical subrootsystem has type $C_3$ realized in $\mathbb{R}^3$ with simple classical roots 
$\alpha_1=(1,-1,0)$, $\alpha_2=(0,1,-1)$, $\alpha_3=(0,0,2)$, simple coroots 
$\alpha_i^\vee=\alpha_i$ for $i\in\{1,2\}$ and $\alpha_3^\vee=(0,0,1)$, and $\varphi=(1,1,0)$. 
We have $2\rho^\vee=(5,3,1)$, $\pair{2\rho^\vee,\varphi}=8$,
$s_\varphi=s_2s_1s_3s_2s_1s_3s_2$ and $s_0 = t_{\varphi} s_{\varphi}$. 
$\quantumBruhatGraph$ has a quantum edge from $s_\varphi$ to $1$ since $\ell(s_\varphi)=7= \pair{2\rho^\vee,\varphi} - 1$.
\item Consider the root system $D_{3+1}^{(2)}$, which is dual to $C_3^{(1)}$.
The classical subrootsystem is type $B_3$ realized in $\mathbb{R}^3$ by $\alpha_1=(1,-1,0)$, $\alpha_2=(0,1,-1)$, and $\alpha_3=(0,0,1)$, $\alpha_i^\vee=\alpha_i$ for $i=1,2$, $\alpha_3^\vee=(0,0,2)$, $2\rho^\vee = (6,4,2)$, $\varphi=(1,0,0)$, $\pair{2\rho^\vee,\varphi}=6$, and $s_\varphi=s_1s_2s_3s_2s_1$. $\quantumBruhatGraph$ has a quantum edge from $s_\varphi$ to $1$ since $\ell(s_\varphi)=5=\pair{2\rho^\vee,\varphi}-1$.
\end{enumerate}
\end{example}


\begin{thebibliography}{1}

\bibitem{comb_of_cox}
Anders Bj{\"{o}}rner and Francesco Brenti.
Combinatorics of Coxeter groups. Graduate Texts in Mathematics, 231. 
Springer, New York, 2005. 

\bibitem{bfp}
Francesco Brenti, Sergey Fomin, and Alexander Postnikov.
Mixed Bruhat operators and Yang-Baxter equations for Weyl groups.
Internat. Math. Res. Notices 1999, no. 8, 419–441. 

\bibitem{deodhar}
Vinay Deodhar.
Some characterizations of Bruhat ordering on a Coxeter group and determination of the relative M\"obius function. Invent. Math. 39 (1977), no. 2, 187–198. 

\bibitem{dyer}
M.~J. Dyer.
Hecke algebras and shellings of Bruhat intervals. Compositio Math. 89 (1993), no. 1, 91–115. 

\bibitem{kumar_email}
Shrawan Kumar.
\newblock Personal correspondence, 2019.

\bibitem{lam_shimo}
Thomas {Lam} and Mark {Shimozono}.
Quantum cohomology of G/P and homology of affine Grassmannian.
Acta Math. 204 (2010), no. 1, 49–90. 

\bibitem{milicevic}
Elizabeth {Mili{\'c}evi{\'c}}.
Maximal Newton polygons via the quantum Bruhat graph. 24th International Conference on Formal Power Series and Algebraic Combinatorics (FPSAC 2012), 899–910,
Discrete Math. Theor. Comput. Sci. Proc., AR, Assoc. Discrete Math. Theor. Comput. Sci., Nancy, 2012. 

\bibitem{postnikov}
Alexander {Postnikov}.
Quantum Bruhat graph and Schubert polynomials. 
Proc. Amer. Math. Soc. 133 (2005), no. 3, 699–709. 

\bibitem{welch}
Amanda Welch.
Double Affine Bruhat Order. 
PhD thesis, Virginia Tech, May 2019.

\end{thebibliography}
\end{document}